\documentclass{article}
\usepackage[english]{babel}
\input xy
\xyoption{all}
\usepackage{amssymb,amsmath,amsthm}

\newcommand{\R}{\mathbb{R}}
\newcommand{\Z}{\mathbb{Z}}

\newcommand{\re}{\mathrm{Re}\,}
\newcommand{\im}{\mathrm{Im}\,}
\newcommand{\lie}[1]{\mathfrak{#1}}     
\newcommand{\Lie}{\mathcal{L}}          
\newcommand{\h}{\mathbb{H}}
\newcommand{\C}{\mathbb{C}}
\newcommand{\CP}{\mathbb{CP}}

\newcommand{\hook}{\lrcorner\,}
\newcommand{\LieG}[1]{\mathrm{#1}}      

\newcommand{\Spin}{\mathrm{Spin}}
\newcommand{\SU}{\mathrm{SU}}
\newcommand{\SO}{\mathrm{SO}}
\newcommand{\Cl}{\mathrm{Cl}}
\newcommand{\Gtwo}{\mathrm{G}_2}
\newcommand{\so}{\mathfrak{so}}
\newcommand{\su}{\mathfrak{su}}

\newcommand{\GL}{\mathrm{GL}}

\newcommand{\dfn}[1]{\emph{#1}}
\newcommand{\id}{\mathrm{Id}}   

\DeclareMathOperator{\Sym}{Sym}

\DeclareMathOperator{\Ad}{Ad}
\DeclareMathOperator{\diag}{diag}

\theoremstyle{plain}
\newtheorem{proposition}{Proposition}
\newtheorem{theorem}[proposition]{Theorem}
\newtheorem{lemma}[proposition]{Lemma}
\newtheorem{corollary}[proposition]{Corollary}
\theoremstyle{definition}

\theoremstyle{remark}
\newtheorem*{remark}{Remark}

\newcommand{\Span}[1]{\operatorname{Span}\{#1\}}
\newcommand{\SpanC}[1]{\operatorname{Span}_\C\{#1\}}

\renewcommand{\Re}{\operatorname{Re}}

\def\cocoa{{\hbox{\rm C\kern-.13em o\kern-.07em C\kern-.13em o\kern-.15em A}}}
\begin{document}
\title{Calabi-Yau cones from contact reduction}
\author{D. Conti and A. Fino}
\date{February 8, 2010}
\maketitle
\begin{abstract}
We consider a generalization of Einstein-Sasaki manifolds, which we characterize in terms both of spinors and differential forms, that in the real analytic case corresponds to contact manifolds whose symplectic cone is Calabi-Yau. We construct solvable examples in seven dimensions. Then, we consider circle actions that preserve the structure, and determine conditions for the contact reduction to carry an induced structure of the same type. We apply this construction to obtain a new hypo-contact structure on $S^2\times\ T^3$.
\end{abstract}
\vskip5pt\centerline{\small\textbf{MSC classification}: 53C25; 53D20, 53C30.}

\centerline{\small\textbf{Keywords}: Contact reduction, generalized Killing spinor, Calabi-Yau.}
\sloppy
\section*{Introduction}
Einstein-Sasaki manifolds can be characterized in terms of real Killing spinors \cite{FriedrichKath}; infinite families of explicit examples of Einstein-Sasaki manifolds have been constructed in \cite{GauntlettMartelliSparksWaldram,CveticEtAl:NewEinsteinSasaki}. These examples are in fact toric  contact manifolds (see \cite{Lerman:ContactToricManifolds}). More generally, one can consider the geometry associated to a generalized Killing spinor \cite{BarGauduchonMoroianu}. Both Killing and generalized Killing spinors can be associated to other geometries (see \cite{Bar}), but we shall only consider generalized Killing spinors associated to $\SU(n)$\nobreakdash-structures on manifolds of dimension $2n+1$. These structures  correspond ideally to restrictions of  Calabi-Yau structures to a hypersurface. In fact, we can think of an $\SU(n)$-structure as given by a real one-form $\alpha$, a real two-form $F$ and a complex $n$-form $\Omega$, and the associated spinor is generalized Killing if and only if
\[dF=0,\quad d(\alpha\wedge\Omega)=0\;.\]
Every hypersurface $M$ in a Calabi-Yau manifold of real dimension $2n+2$ has an $\SU(n)$-structure of this type, where $F$ and $\alpha\wedge\Omega$ are the restriction of the K\"ahler form and complex volume. In the real analytic category, the converse also holds  (see \cite{ContiSalamon} for the five-dimensional case, and Proposition~\ref{prop:GKEqualsHypo} for the general case).

Realizing $M$ explicitly as a hypersurface in a Calabi-Yau manifold is generally a matter of solving certain evolution equations; in the  Einstein-Sasaki case, however, the Calabi-Yau structure is simply the induced conical \mbox{$\SU(n+1)$}\nobreakdash-structure on $M\times\R_+$. Thus, the Calabi-Yau manifold $M\times\R_+$ is both the Riemannian  and  symplectic cone over $M$. In this paper we consider the intermediate case in which we have a generalized Killing spinor, and the $\SU(n)$-structure on $M$ is a contact metric structure, thereby inducing a conical symplectic form on $M\times\R_+$; by \cite{FriedrichKim:TheEinsteinDiracEquation}, $\alpha$-Einstein-Sasaki manifolds belong to this class. In this context, \emph{contact} means that $d\alpha=-2F$. If $(M,\alpha)$ is a real analytic contact manifold with a real analytic $\SU(n)$-structure as above, it turns out that the symplectic cone $M\times\R_+$ has a compatible Calabi-Yau metric; in other words, the Calabi-Yau metric whose existence is guaranteed by the above-mentioned  embedding result has conical K\"ahler form, although the metric itself may not be conical (Theorem~\ref{thm:HypoContact}).

Not many examples of $\SU(n)$-structures of this type are known, excepting those that are actually $\alpha$-Einstein-Sasaki.  There are  solvable non-compact examples in five dimensions  \cite{DeAndres:HypoContact}, an example in seven dimensions related to generalized $\Gtwo$-structures \cite{FinoTomassini}, and a two-parameter family of examples in the sphere bundle in $T\CP^2$ \cite{Conti:InvariantForms}. Notice however that examples of $\alpha$-Einstein-Sasaki manifolds which are not Einstein are known, the simplest being the Heisenberg group in arbitrary odd dimension \cite{TomassiniVezzoni}; more complicated $\alpha$-Einstein-Sasaki structures were constructed in \cite{BoyerGalickiMatzeu}.

\medskip
The $7$-dimensional example  in \cite{FinoTomassini} is constructed as a compact quotient of the Lie group $\SU(2) \ltimes_{\varphi} \R^4$, where $\varphi$ is given by quaternionic multiplication. In Section~\ref{sec:ExamplesNoQuotient} we show that this example is unique in the class of semidirect products $H \ltimes V$ of dimension higher than three, with $H$ compact Lie group and $V$ representation of $H$.  Motivated by this uniqueness,  we consider semidirect products $H\ltimes V$ with $H$ solvable, obtaining new seven-dimensional examples of contact $\SU(3)$-structures associated to a generalized Killing spinor. In fact, we classify the three-dimensional solvable Lie groups $H$ which give rise to structures of this type on a semidirect product $H\ltimes \R^4$ (Theorem~\ref{thm:Classification}).

Our main result is aimed at constructing more examples by using the contact reduction. We start with a manifold $M$ of dimension $2n+1$, a generalized Killing spinor $\psi$ and an $\SU(n)$\nobreakdash-structure which is a contact metric structure. Moreover, we assume that $S^1$ acts on $M$ preserving the $\SU(n)$-structure. Under certain conditions on this action,  one can define the contact reduction $M//S^1$, which is a contact manifold of dimension $2n-1$ \cite{Geiges:Constructions}, and inherits a Riemannian metric and spinor as well. In Theorem~\ref{thm:HypoToHypo} we prove a necessary and sufficient condition for this induced spinor to be generalized Killing; more precisely, this condition involves the derivative of the norm of the fundamental vector field, and it is satisfied by circle actions with constant orbit length. The proof depends on a study of ``basic'' spinors in the context of Riemannian submersions, which we carry out in Section~\ref{sec:Submersions}.

Whilst serving as motivation for the study of these structures, the existence of a Calabi-Yau cone is not used directly; in particular, our main result holds in the smooth category. However, if  Calabi-Yau cones do exist, then the cone over $M//S^1$ is the symplectic reduction of the cone over $M$. It is well known that the symplectic reduction of a K\"ahler manifold is again K\"ahler \cite{HitchinKarlhedeLindstromRovcek}, but Ricci-flatness is not preserved in general. Thus, our result gives in particular sufficient conditions for the Calabi-Yau condition to be preserved under K\"ahler reduction. The odd-dimensional analogue of K\"ahler reduction was considered in \cite{GrantcharovOrnea}, where it was shown that the Sasaki condition is preserved under contact reduction. Neither is the Einstein condition preserved in this situation; however, our result gives sufficient conditions for the Sasaki reduction of a Einstein-Sasaki manifold to be $\alpha$-Einstein (Corollary~\ref{cor:alphaEinsteinSasaki}).

As an application of our theorem, we obtain a new example on the compact manifold $S^2\times T^3$.

\section{Generalized Killing spinors}
\label{sec:GKSpinors}
Let $M$ be an oriented spin manifold of dimension $n$; the choice of a Riemannian metric $g$ on $M$ determines the principal bundle $P_{\SO}$ of oriented orthonormal frames. More generally, if $G$ is a subgroup of $\GL(n,\R)$, a $G$-structure $P_G$ on $M$ is a reduction to $G$ of the bundle of frames. In this section we shall consider two types of $G$-structures associated to a spinor, corresponding to the subgroups $\SU(n)\subset\GL(2n,\R)$ and $\SU(n)\subset\GL(2n+1,\R)$. The first type corresponds to Calabi-Yau geometry; more precisely, an $\SU(n)$-structure $P_\SU$ on a manifold of dimension $2n$ is called \dfn{Calabi-Yau} if it admits a torsion-free connection. In particular, a Calabi-Yau structure determines a Riemannian metric with holonomy contained in $\SU(n)$.

The Lie group $\Spin(n)$ acts on $P_{\SO}$ on the right via the 2:1 homomorphism
\begin{equation}
\label{eqn:Ad}
\Ad\colon \Spin(n)\to\SO(n)\;.
\end{equation}
A \dfn{spin structure} on $M$ is a principal bundle over $M$ with fibre $\Spin(n)$ and an equivariant 2:1 map
\begin{equation*}
\label{eqn:BundleAd}
\Ad\colon P_{\Spin}\to P_{\SO}\;.
\end{equation*}
Let $\Sigma=\Sigma^n$ be a complex irreducible representation of the Clifford algebra \mbox{$\Cl(n)\supset\Spin(n)$};
a \dfn{spinor} on $M$ is a section of the associated bundle
\[P_{\Spin}\times_{\Spin(n)}\Sigma\;.\]

We shall need explicit formulae for the representation $\Sigma$. Let $v_0,v_1$ be a basis of $\C^2$ with
$v_0=\left(\begin{smallmatrix}1\\0\end{smallmatrix}\right)$ and
$v_1=\left(\begin{smallmatrix}0\\1\end{smallmatrix}\right)$; a basis $u_0,\dots,u_{2^n-1}$ of $\Sigma=(\C^2)^{\otimes
n}$ is given by
\begin{align*}
u_k&=v_{a_{n-1}}\otimes\dots\otimes v_{a_0}\;,&\text{ where }\quad k=\sum_{0\leq r< n} a_r2^r\;.
\end{align*}
We think of the $a_r$ as elements of $\Z/2\Z$. The Clifford algebra $\Cl(2n)$ acts irreducibly on $\Sigma$ by
\begin{align*}
e_{2j}\cdot u_k&=-(-1)^{a_{j-1}+\dots+a_0} v_{a_{n-1}}\otimes\dots\otimes v_{a_{j}}\otimes v_{1-a_{j-1}}\otimes
v_{a_{j-2}}\otimes\dots\otimes v_{a_0}\\
e_{2j-1}\cdot u_k&=i(-1)^{a_{j-2}+\dots+a_0} v_{a_{n-1}}\otimes\dots\otimes v_{a_{j}}\otimes v_{1-a_{j-1}}\otimes
v_{a_{j-2}}\otimes\dots\otimes v_{a_0}
\end{align*}
where $j$ ranges between $1$ and $n$.
The representation $\Sigma$ is also an irreducible representation of $\Cl(2n+1)$; more precisely, denoting temporarily by $\odot$ the Clifford multiplication of $\Cl(2n)$, we set
\begin{align*}
e_{2n+1}\cdot u_k&=(-1)^{a_{n-1}+\dots+a_0+n}i\,u_k,\\
 e_j\cdot u_k&=-e_{2n+1}\cdot e_j\odot u_k \;,& j=1,\dotsc,2n.
\end{align*}
There is a choice of sign involved in this definition; in this paper we shall only consider this representation of $\Cl(2n+1)$, which is characterized by the fact that the volume element $e_1\dotsm e_{2n+1}$ acts as $i^{n+1}$.

Restricting the  action to the groups $\Spin(2n)$, $\Spin(2n+1)$, one finds that the stabilizer of $u_0$ is $\SU(n)$ in both cases.
However,  $\Sigma$ is only irreducible as a $\Spin(2n+1)$-module, whilst as a $\Spin(2n)$-module it splits into the two components
 \begin{align*}
 \Sigma_+&=\Span{ u_k\mid a_{n-1}+\dots+a_0=0 \mod 2}\;,\\
 \Sigma_-&=\Span{u_k\mid a_{n-1}+\dots+a_0=1 \mod 2}\;.
 \end{align*}

As a consequence of the above, Calabi-Yau manifolds, namely Riemannian manifolds of dimension $2n$ with holonomy $\SU(n)$, can be characterized by the existence of a parallel spinor. Indeed, the orbit $\Spin(2n)u_0$ defines a bundle
\begin{equation}
\label{eqn:OneOrbit}
P_{\Spin}\times_{\Spin(2n)}\Spin(2n)u_0\subset P_{\Spin}\times_{\Spin(2n)}\Sigma\;,
\end{equation}
and a section $\psi$ of this bundle defines a principal bundle
\[P_\SU = \{u\in P_\Spin \mid [u,u_0]=\psi\},\]
which has the stabilizer $\SU(n)$ as its fibre.
Since the kernel of \eqref{eqn:Ad} is not contained in $\SU(n)$, the principal bundle $P_\SU\cong \Ad(P_\SU)$ is in fact an $\SU(n)$-structure. Hence, the covariant derivative of $\psi$ with respect to the Levi Civita connection can be identified with the intrinsic torsion of $P_\SU$. In particular, if  $\psi$ is parallel, one obtains a Calabi-Yau structure on $M$.

There are other, weaker conditions that it is natural to impose on a spinor. The authors of \cite{BarGauduchonMoroianu} propose the following: a \dfn{generalized Killing spinor} $M$ is a spinor  $\psi$ satisfying
\begin{equation}
 \label{eqn:GKSpinor}
\nabla_X\psi =\frac12Q(X)\cdot\psi\;,\quad Q\in\Gamma(\Sym(TM))\;.
\end{equation}
Here, as throughout the paper, $\nabla$ denotes the Levi-Civita connection. Generalized Killing spinors arise naturally, by restriction, on oriented hypersurfaces in manifolds with a parallel spinor; in this setting, the tensor $Q$ corresponds to the Weingarten tensor. There are also partial results in the converse direction (see \cite{Morel,BarGauduchonMoroianu,ContiSalamon}).

In particular, consider the case of a hypersurface $M$ inside a manifold with holonomy $\SU(n+1)$. The restriction of the parallel spinor gives a section of \eqref{eqn:OneOrbit} on $M$, thus defining an $\SU(n)$-structure $P_\SU$ on the $2n+1$-dimensional manifold $M$. To this $\SU(n)$-structure one can associate forms $\alpha$, $F$ and $\Omega$ such that
\begin{equation}
\label{eqn:StandardFrame}
\begin{gathered}
\alpha=e^{2n+1},\quad F=e^{12}+\dotsb+e^{2n-1,2n},\\ \Omega=(e^1+ie^2)\wedge\dotsb\wedge (e^{2n-1}+ie^{2n}),
\end{gathered}
\end{equation}
where $e^1,\dotsc, e^{2n+1}$ is the coframe associated to to any local section of $P_\SU$, and $e^{jk}$ is short for $e^j\wedge e^k$. One can also read off the forms $F$ and $\Omega\wedge\alpha$  as pull-backs of the K\"ahler form and the complex volume form on the Calabi-Yau manifold. We shall say that $P_\SU$ is the restriction of the Calabi-Yau structure to the hypersurface.

Regardless of whether the $\SU(n)$-structure on $M$ arises from restricting a Calabi-Yau structure, assuming \eqref{eqn:GKSpinor} holds, one can relate the covariant derivative of the forms $\alpha$, $F$ and $\Omega$ to the tensor $Q$. This will enable us to rewrite the generalized Killing spinor equation \eqref{eqn:GKSpinor} in terms of differential forms in the next section.
\begin{lemma}
\label{lemma:NablaOfForms}
Let $M$ be a Riemannian spin manifold of dimension $2n+1$, and let $\psi$ be a section of \eqref{eqn:OneOrbit} satisfying \eqref{eqn:GKSpinor}. If $\alpha$, $F$ and $\Omega$ are the forms associated to the $\SU(n)$-structure defined by $\psi$, then
\begin{gather*}
\nabla_X\alpha =  (-1)^{n}Q(X)\hook F\\
\nabla_X F=(-1)^n\alpha\wedge Q(X)\\
\nabla_X \Omega=(-1)^{n+1}i\alpha\wedge (Q(X)\hook\Omega)+i (-1)^n Q(X,e_{2n+1})\Omega
\end{gather*}
where we have used the identification $TM\cong T^*M$ given by the metric.
\end{lemma}
\begin{proof}
We work on the principal bundle $P_\SU\cong\Ad(P_\SU)$. Let $\omega$ be the restriction to $P_\SU$ of the Levi-Civita connection form. According to the orthogonal decomposition
\[\so(2n+1)=\su(n)\oplus\langle J\rangle\oplus \lie{u}(n)^\perp\;,\]
where
\[J=e_{21}-e_{12}+\dotsb+e_{2n,2n-1}-e_{2n-1,2n}\;,\]
we can decompose $\omega$  as $\omega_\su + kJ+\omega^\perp$. Here, $e_{hk}$ is the square matrix of order $2n+1$ with entries $a_{ij}$, all equal to zero except $a_{hk}=1$.
Then
\[\nabla_X \psi = \omega(X)u_0=k(X)J u_0+\omega^\perp(X)u_0=-\frac12ni\,k(X)u_0-\frac12 (\omega^\perp)_{ij} e_i\cdot e_j\cdot u_0,\]
where \[\omega^\perp=\sum_{1\leq i<j\leq 2n+1} (\omega^\perp)_{ij}(e_{ij}-e_{ji})\;.\]
Now, by \eqref{eqn:GKSpinor},
\[ni\,k(X) \,u_0-(\omega^\perp)_{ij}(X) e_i\cdot e_j\cdot u_0=Q(X)\cdot u_0;\]
looking at the definition of the Clifford action, from which it follows  in particular that for $1\leq j\leq n$
\begin{align*}
e_{2j}\cdot e_{2n+1}\cdot u_0&=-e_{2n+1}\cdot e_{2j}\cdot u_0=(-1)^{n+1}e_{2j-1}\cdot u_0,\\
e_{2j-1}\cdot e_{2n+1}\cdot u_0 &= -e_{2n+1}\cdot e_{2j-1}\cdot u_0 =(-1)^{n}e_{2j}\cdot u_0,
\end{align*}
one concludes that
\begin{gather*}
(\omega^\perp)_{ij}=0,\quad j<2n+1,\quad  (\omega^\perp)_{i,2n+1}(X)=(-1)^{n+1}Q(X,J(e_i)),\\
k(X)=-\frac1n (-1)^n Q(X,e_{2n+1}).
\end{gather*}
Now observe that
\[J\alpha=0\,\quad J F=0,\quad J\Omega = -ni\Omega;\]
hence,
\[\nabla_X\alpha = \omega(X)\alpha =\omega^\perp(X)\alpha = (-1)^{n}J(Q(X))\;.\]
The rest of the statement is proved in the same way.
\end{proof}
\begin{remark}
The forms $\alpha$, $ F$ and $\Omega$ determine the $\SU(n)$-structure $P_\SU$. Therefore, one can express the intrinsic torsion of $P_\SU$ in terms of $\nabla\alpha$, $\nabla F$ and $\nabla\Omega$; however, the intrinsic torsion turns out to be entirely determined by $d\alpha$, $d F$ and $d\Omega$ (see \cite{Conti:Embedding}).
\end{remark}
\section{Calabi-Yau cones}
\label{sec:Cones}
In this section we restrict to the real analytic category, and give a twofold chacterization of contact $\SU(n)$-structures associated to a generalized Killing spinor, in terms of differential forms and Calabi-Yau cones.

First, we need to consider a broader class of $\SU(n)$-structures.
\begin{proposition}
\label{prop:GKEqualsHypo}
Let $M$ be a real analytic manifold of dimension $2n+1$ with a real analytic $\SU(n)$-structure $P_\SU$. The following are equivalent:
\begin{enumerate}
 \item The section of the vector bundle  \eqref{eqn:OneOrbit} associated to $P_\SU$ is a generalized Killing spinor.
\item The differential forms $\alpha$, $ F$ and $\Omega$ associated to $P_\SU$ satisfy
\[d F=0,\quad d(\alpha\wedge\Omega)=0.\]
\item A neighbourhood of $M\times\{0\}$ in $M\times\R$ has a Calabi-Yau structure which restricts to $P_\SU$.
\end{enumerate}
\end{proposition}
\begin{proof}
Assume that (i) holds, and let $e_i$ be the local orthonormal frame associated to a section of $P_\SU$. Since $Q$ is symmetric, it follows that
\[\sum_i e^i\wedge Q(e_i)=0.\]
By Lemma~\ref{lemma:NablaOfForms}, and using the fact that the Levi-Civita connection is torsion-free,
\[d F=\sum_i e^i\wedge \nabla_{e_i}  F=(-1)^n\sum_i e^i\wedge \alpha\wedge Q(e_i)=0.\]
Likewise,
\[d(\alpha\wedge\Omega)=(-1)^n \sum_i e^i\wedge \left((Q(e_i)\hook F) +i Q(e_i,e_{2n+1}) \alpha\right)\wedge\Omega\;;\]
since for any vector field $Y$ one has
\[\bigl(Y\hook F-iY+i\alpha(Y)\alpha\bigr)\wedge\Omega=0,\]
we find that
\[d(\alpha\wedge\Omega)=(-1)^n \sum_{i} e^i\wedge iQ(e_i)\wedge\Omega=0,\]
meaning that (ii) holds.

The fact that (ii) implies (iii) follows from the theory of exterior differential systems. Reference \cite{ContiSalamon} contains a detailed proof of the five-dimensional case, which can be generalized to arbitrary dimension because the exterior differential system associated  to the group $\SU(n)$ is involutive \cite{Bryant:Calibrated,Conti:Embedding}.

Finally, (iii) implies (i) because every Calabi-Yau manifold carries a parallel spinor $\psi$, and its restriction to a hypersurface satisfies
\[\nabla_X \iota^*\psi = \frac12 Q(X)\cdot\iota^*(\nabla_X \psi),\]
where $Q$ is the Weingarten tensor (see e.g. \cite{BarGauduchonMoroianu}).
\end{proof}
\begin{remark}
The assumption of real analyticity in Proposition~\ref{prop:GKEqualsHypo} is certainly necessary to prove that (i) or (ii) imply (iii), due to a result of Bryant \cite{Bryant:EmbeddingProperty}. However, the fact that (i) implies (ii) does not require this hypothesis. A five-dimensional version of Proposition~\ref{prop:GKEqualsHypo} is proved in \cite{ContiSalamon}, where the $\SU(2)$-structures defined by a generalized Killing spinor were introduced under the name of hypo structures. In that paper it was also proved, by considering the intrinsic torsion, that (ii) implies (i) without assuming real analyticity. We expect a similar result to hold in arbitrary dimension.
\end{remark}
\begin{remark}
The passage from (ii) to (iii) can be described in terms of evolution equations, in the sense of \cite{Hitchin:StableForms}. Indeed, suppose  there is a one-parameter family $(\alpha(t),F(t),\Omega(t))$,  of $\SU(n)$-structures on $M$, with $t$ ranging in the interval $(a,b)$; then the forms
\begin{equation}
\label{eqn:CYFromOneParameterFamily}
\alpha(t)\wedge dt + F(t),\quad (\alpha(t)+idt)\wedge\Omega(t)
\end{equation}
define a Calabi-Yau structure on $M\times (a,b)$ if and only if (ii) holds for, say, $t=0$ and the evolution equations
\begin{equation}
\label{eqn:EvolutionEquations}
\frac{\partial}{\partial t} F(t) =-d\alpha(t),\quad \frac{\partial}{\partial t} (\alpha(t)\wedge\Omega(t))=id\Omega(t)
\end{equation}
are satisfied. Conversely, if $M$ is compact (iii) implies that one can find a solution of these equations with $(\alpha(0),F(0),\Omega(0))$ corresponding to $P_\SU$. Indeed,  the exponential map enables one to identify a tubular neighbourhood of $M$ in $M\times\R$ with a product $M\times(a,b)$, in such a way that the vector field $\frac{\partial}{\partial t}$ has unit norm and is orthogonal to the hypersurfaces $M\times\{t\}$, with the effect of casting the K\"ahler form and complex volume  in the form \eqref{eqn:CYFromOneParameterFamily}. Thus, Proposition~\ref{prop:GKEqualsHypo} can be viewed as an existence result for solutions of \eqref{eqn:EvolutionEquations}.
\end{remark}

A special situation of Proposition~\ref{prop:GKEqualsHypo} is when $Q$ is a (constant multiple of) the identity. Then the spinor $\psi$ is called a real Killing spinor, and it is the restriction of a parallel spinor on the Riemannian cone over $M$. The general situation has been studied in  \cite{Bar}; in our case, the restriction of the Calabi-Yau structure is an Einstein-Sasaki structure. In particular, this means that  $M$ is a contact metric manifold, with contact form $\alpha$. We define the \dfn{symplectic cone} over $(M,\alpha)$  as the symplectic manifold
\begin{equation}
 \label{eqn:SymplecticCone}
\left(M\times\R_+,-\frac12 d(r^2\alpha)\right).
\end{equation}
If $M$ is Einstein-Sasaki, the symplectic cone is Calabi-Yau with the cone metric $r^2g+dr^2$; it is understood, here and wherever we refer to Calabi-Yau structures on symplectic manifolds,  that  the K\"ahler form coincides with the given symplectic form.

More generally, we  say that an  $\SU(n)$-structure on a manifold $M$ of dimension $2n+1$ is \dfn{contact} if $d\alpha=-2 F$; this means that $\alpha$ is a contact form, and $F$ is the pullback to $M\cong M\times\{1\}$ of the conical  symplectic form \eqref{eqn:SymplecticCone}.
We shall consider a weaker condition than Einstein-Sasaki, corresponding ideally to contact $\SU(n)$-structures $P_\SU$  such that the symplectic cone is Calabi-Yau, but not necessarily with respect to the cone metric.

\begin{theorem}
\label{thm:HypoContact}
Let $M$ be a real analytic manifold of dimension $2n+1$ with a real analytic, contact $\SU(n)$-structure $P_\SU$. The following are equivalent:
\begin{enumerate}
 \item The section of  \eqref{eqn:OneOrbit} associated to $P_\SU$ is a generalized Killing spinor.
\item The differential forms $\alpha$, $ F$ and $\Omega$ associated to $P_\SU$ satisfy
\[d\alpha=-2 F,\quad \alpha\wedge d\Omega=0.\]
\item A neighbourhood of $M\times\{1\}$ in the symplectic cone $M\times\R_+$ has a Calabi-Yau metric which restricts to $P_\SU$.
\end{enumerate}
\end{theorem}
\begin{proof}
The fact that (i) implies (ii) and (iii) implies (i) is a consequence of Proposition~\ref{prop:GKEqualsHypo}.

To see that (ii) implies (iii), one applies Proposition~\ref{prop:GKEqualsHypo}, deducing that a neighbourhood $N$ of $M\times\{0\}$ in $M\times\R$
 has a Calabi-Yau structure restricting to $P_\SU$ on $M\cong M\times\{0\}$. Under the diffeomorphism
\[M\times\R\ni(x,t)\to (x,e^t)\in M\times\R_+,\]
the conical symplectic form reads
\[\omega_0=e^{2t}(\alpha\wedge dt+F).\]
Thus, it suffices to prove that there is a  diffeomorphism of a neighbourhood $N'$ of $M\times\{0\}$ into $N$ that is the identity on $M\times\{0\}$ and pulls back the Calabi-Yau symplectic form into the conical symplectic form.

Let $\alpha$, $F$, $\Omega$ be the forms on $N$  given by the restriction of the Calabi-Yau structure to each hypersurface $N\cap (M\times\{t\})$. By construction, the K\"ahler form is given by
\[\omega_1=\alpha\wedge dt+F.\]
Consider a time-dependent vector field $X_s$ on $N$, and let $\phi_s$
be the flow of $X_s$, which satisfies
\[\phi_0(x,t)=(x,t),\quad \frac{\partial}{\partial s}\phi_s(x,t)=(X_s)_{\phi_s(x,t)}\;.\]
It is a general fact that for every form $\beta$, one has
\begin{equation}
\label{eqn:GeneralFact}
\phi_{s}^* \Lie_{X_s} \beta = \frac{\partial}{\partial s} \phi_{s}^*\beta.
\end{equation}
Up to restricting $N$, we can define a one-parameter family of symplectic forms by
\[(\omega_s)_{(x,t)} = (1-s)e^{2t} \left(\alpha(x,0)\wedge dt+ F(x,0)\right) + s\left(\alpha(x,t)\wedge dt + F(x,t)\right),\]
thus interpolating between the conical symplectic form $\omega_0$ and the K\"ahler form $\omega_1$.
We shall determine $X_s$ in such a way that $\phi_s\colon N'\to N$ is well defined for $0\leq s\leq 1$, and
\begin{equation}
\label{eqn:omegas}
\omega_0=\phi_s^*\omega_s.
\end{equation}
Equality certainly holds for $s=0$; taking the derivative with respect to $s$ and applying \eqref{eqn:GeneralFact}, we get
\[0=\phi_{s}^* \left(\Lie_{X_s}\omega_s + \omega_1-\omega_0\right)=\phi_{s}^* \left(d(X_s\hook\omega_s) +\omega_1-\omega_0\right).\]
The two-form $\omega_1-\omega_0$ is cohomologically trivial on $N$, because it vanishes on  $M\times\{0\}$, that we may assume to be a deformation retract of $N$. It follows that $\omega_1-\omega_0=d\beta$ for some $1$-form $\beta$. Thus, it suffices to require
\[X_s \hook\omega_s = -\beta,\]
which determines $X_s$ because the $\omega_s$ are symplectic forms, to ensure that \eqref{eqn:omegas} holds. With this definition, $X_s$ vanishes on $M\times\{0\}$. Hence, the flow $\phi_s$ is well defined for all $s$ (indeed, constant) at $t=0$, and up to restricting $N'$, we can assume that $\phi_s$ is well defined for $0\leq s\leq 1$.
Now set $\Phi=\phi_1\colon N'\to N$; this is a diffeomorphism that pulls back the K\"ahler form to the conical symplectic form. Since $\Phi(x,0)=(x,0)$,
we have the following diagram
\[
 \xymatrix{(N',\omega_0)\ar[r]^\Phi & (N,\omega_1)\\
		M\times\{0\}\ar[r]^{\id}\ar[u] & M\times\{0\}\ar[u]
}
\]
where the vertical arrows are inclusions. Since the diagram commutes, the pullback under $\Phi^*$ of the Calabi-Yau structure on $M\times\R$ restricts to $P_\SU$ on $M\times\{0\}$.
\end{proof}
\begin{remark}
Again, the assumption of real analyticity is only essential to prove the implication (ii) $\implies$ (iii). We shall not need this hypothesis in the rest of the paper.
\end{remark}
\begin{remark}
 In the proof of Theorem~\ref{thm:HypoContact}, the Calabi-Yau structure on $N$ can be described in terms of the evolution equations \eqref{eqn:EvolutionEquations}, by requiring that each $M\times\{t\}$ be orthogonal to the unit vector field $\frac{\partial}{\partial t}$. However, the map $\Phi$ does not preserve this description. Thus, when working on $N'$, the one-parameter family of $\SU(n)$-structures induced by the inclusions $M\times\{t\}\subset N'$ will not satisfy, in general, the  evolution equations. From the side of $N$, this means that one should not expect ``conical'' evolution.
\end{remark}

Hypo-contact manifolds, namely five-dimensional manifolds with $\SU(2)$\nobreakdash-structures satisfying Theorem~\ref{thm:HypoContact}, have been studied in \cite{DeAndres:HypoContact}, which contains a classification of  solvable Lie groups with invariant hypo-contact structures.

\section{Seven-dimensional semidirect products}
\label{sec:ExamplesNoQuotient}
In this section we give new examples of $\SU(3)$\nobreakdash-structures in seven dimensions  satisfying Theorem~\ref{thm:HypoContact}.
More precisely, we consider semidirect products $H\ltimes V$, with $H$ a Lie group, and $V$ a representation of $H$, generalizing the example $\SU(2)\ltimes\R^4$ of \cite{FinoTomassini} (also reviewed in Section~\ref{sec:Examples}). We show that this example is unique among those with $H$ compact and connected, at least when the overall dimension is higher than three. We then classify  the solvable $3$-dimensional Lie groups $H$ such that some semidirect product $H\ltimes \R^4$ admits a left-invariant contact $\SU(3)$-structure whose associated spinor is generalized Killing.

\begin{proposition}
\label{prop:CompactUnique}
Let $H$ be a compact connected Lie group, and $\varphi\colon H\to V$ a representation of $V$. Then the semidirect product $H\ltimes_\varphi V$ has a left-invariant contact structure if and only if $H\ltimes_\varphi V$ is either $\SU(2)\ltimes_\varphi\h$ or $\LieG{U}(1)\ltimes\C$.
\end{proposition}
\begin{proof}
By definition, $H\ltimes_\varphi V$ is the product $H\times V$ with multiplication law given by
\[(h,v)(h',v')=(hh',\varphi(h)v'+v).\]
Now let $\alpha$ be a left-invariant one-form. We can write $\alpha_{(e,0)}=(\alpha_H)_e+\alpha_V$, where $\alpha_H$ is a left-invariant one-form on $H$ and $\alpha_V$ is in $V^*$. Then
\begin{equation}
 \label{eqn:thecontactform}
\alpha_{(h,v)}=L_{(h,v)}^*(\alpha_H+\alpha_V) =  L_{h}^*\alpha_H +\varphi(h^{-1})\alpha_V.
\end{equation}
Hence if $e_1,\dotsc,e_n$ is a basis of $\lie{h}$ and $e^1,\dotsc, e^n$ the dual basis of 1-forms,
at $h=e$, we find
\[(d\alpha)_{(e,v)} = (d\alpha_H)_e -  \sum_{i=1}^n e^i\wedge (\varphi_{*e}e_i) \alpha_V.\]
If $\alpha$ is contact, then
\[\dim \{(\varphi_{*e}e_i) \alpha_V\}=\dim \left\{X\hook (d\alpha)_{(e,v)}\mid X\in V\right\}\geq\dim V- 1.\]
Hence, the orbit of $\alpha_V$ in $V^*$ has at most codimension $1$, and since $H$ is compact, this means that $H$ acts transitively on the sphere in $V^*$.

Now let $K$ be the stabilizer of $\alpha_V$ in $V^*$; suppose by contradiction that $K$ is not discrete, and take $k$ in $K$. Then \eqref{eqn:thecontactform} gives
\[\alpha_{(k,v)}= L_{k}^*\alpha_H +\alpha_V.\]
Restricting $\alpha$ to $K\ltimes_\varphi V$, we find
\[d\alpha|_{K\ltimes_\varphi V}= (d\alpha_H)|_K.\]
Now consider the element
\[\beta=(d\alpha|_{K\ltimes_\varphi V})_e\in \Lambda^2(\lie{k}\oplus V)^*;\]
by above, $\beta$ is really contained in $\Lambda^2\lie{k}^*$. Moreover, since $K$ is compact, the exact form $d\alpha_K$ cannot be a symplectic form on $K$, and so there is a subspace $W\subset \lie{k}$ with
\[\beta|_W=0,\quad \dim W>\frac12\dim K.\]
In particular, the restriction of $d\alpha_e$ to $V\oplus W\subset T_{(e,0)}(H\ltimes V)$ is zero.
Since $H/K$ is the sphere in $V$, we have
\[\dim H=\dim K+\dim V-1,\]
and therefore
\[\dim V+\dim W>\dim V+\frac12\dim K =\frac12(\dim V+ \dim H)+\frac12.\]
Thus, $\alpha$ is not contact.

We have shown that, if $\alpha$ is contact, then $H$ acts transitively on the sphere in $V$ with discrete stabilizer $K$. If $V$ has dimension two, this implies trivially that $H=\LieG{U}(1)$. If $V$ has dimension $n>2$, the exact homotopy sequence
\[0=\pi_1(S^{n-1})\to\pi_0(H)\to\pi_0(K)\to 0\]
implies that $K$ is connected, and therefore trivial. Then $H$ acts on $S^{n-1}$ both transitively and freely, giving $H=\SU(2)$ and $V=\h$.
\end{proof}

Now let $H$ be a solvable $3$-dimensional Lie group, and let $e^1,e^4,e^6$ be an invariant basis of one-forms, where the indices have been chosen for compatibility with \eqref{eqn:StandardFrame} (see Lemma~\ref{lemma:SemidirectStructure} below). Up to a change of basis (see \cite{DeGraaf}),  the structure equations of $H$ are given by exactly one of the following:
\begin{align}
\label{eqn:h:abelian} de^1&=0,& de^4&=0,& de^6&=0;\\
\label{eqn:h:L4} de^1&=0,& de^4&=\pm e^{16},& de^6&=e^{14};\\
\label{eqn:h:L2} de^1&=0,& de^4&=e^{14},& de^6&=e^{16};\\
\label{eqn:h:L3:A0} de^1&=0,& de^4&=0,& de^6&=e^{16};\\
\label{eqn:h:nilpotent} de^1&=0,& de^4&=0,& de^6&=e^{14};\\
\label{eqn:h:L3} de^1&=0,& de^4&=Ae^{14},& de^6&=e^{14}+e^{16}.
\end{align}
In \eqref{eqn:h:L3}, $A$ is a non-zero real constant; the case $A=0$ corresponds to \eqref{eqn:h:L3:A0} under a change of basis.

We can now state the main result of this section.
\begin{theorem}
\label{thm:Classification}
Let $H$ be a solvable $3$-dimensional Lie group. Then there exists a semidirect product $H\ltimes\R^4$ admitting a left-invariant contact $\SU(3)$-structure whose associated spinor is generalized Killing if and only if $H$ has structure equations \textup{(\ref{eqn:h:abelian}\nobreakdash--\ref{eqn:h:L3:A0})}.
\end{theorem}
The rest of this section consists in the proof of Theorem~\ref{thm:Classification}. We shall start by giving a ``structure lemma''. As a preliminary observation, notice that ``rotating'' $\Omega$, i.e. multiplying it by a constant $e^{i\theta}$, preserves the condition of Theorem~\ref{thm:HypoContact}.
\begin{lemma}
\label{lemma:SemidirectStructure}
Let $H$ be a three-dimensional Lie group and $V$ a four-dimensional representation of $H$. Suppose the semidirect product $G=H\ltimes V$ has a left-invariant $\SU(3)$-structure $(\alpha,F,\Omega)$, with
\[d\alpha=-2F,\quad d\Omega\wedge\alpha=0.\]
Then, up to rotating $\Omega$, we can choose a left-invariant basis $e^1,\dotsc,e^7$ of one-forms on $G$ satisfying \eqref{eqn:StandardFrame}, with
\begin{gather}
\label{eqn:de2357}
de^2,de^3,de^5\in\Span{e^2,e^3,e^5,e^7}\wedge\Span{e^1,e^4,e^6},
\end{gather}
and $e^1,e^4,e^6$ invariant extensions of forms on $H\ltimes\{0\}\subset G$.
\end{lemma}
\begin{proof}
By invariance, we can work on the Lie algebra $\lie{g}$ of $G$. As a vector space $\lie{g}=\lie{h}\oplus V$, and the Lie bracket satisfies
\begin{equation}
\label{eqn:hPlusV}
[\lie{h},\lie{h}]\subset\lie{h},\quad [\lie{h},V]\subset V, \quad [V,V]=0.
\end{equation}
We first show that the characteristic vector field is in $V$. Consider the linear map
\[\lie{g}\ni X\xrightarrow{\phi} X\hook F\in\lie{g}^*\;;\]
by construction, the kernel of $\phi$ is spanned by the characteristic vector field $e_7$. Since the restriction of  the exact form $F$ to the abelian Lie algebra $V$ is zero, $\phi(V)$ is contained in $\lie{h}^*$; by a dimension count, this implies that $\ker \phi\subset V$, and so the characteristic vector field is in $V$.

Now consider the restriction of $\Omega$ to $V$. This is a complex $3$-form satisfying
\[e_7\hook (\Omega|_V)=0.\]
Since $V$ has dimension four, this means that $\re\Omega|_V$ and $\im\Omega|_V$ are linearly dependent, and one can multiply $\Omega$ by some $e^{i\theta}$, obtaining $\re\Omega|_V=0$. In other words, $(e_7)^\perp\cap V$ is special Lagrangian in $(e_7)^\perp$. Now recall that the structure group $\SU(3)$ acts transitively on special Lagrangian subspaces of $\C^3$; hence,  we can complete $e_7$ to a basis $e_1,\dotsc,e_7$ of $\lie{g}$, consistent with \eqref{eqn:StandardFrame}, such that
\[(e_7)^\perp\cap V=\Span{e_2,e_3,e_5}.\]
Accordingly, $\lie{h}^*=\Span{e^1,e^4,e^6}$, so
by \eqref{eqn:hPlusV}
\begin{equation}
 \label{eqn:hPlusV2}
d(\lie{g}^*)\subset \Span{e^{1},e^{4},e^{6}}\wedge \lie{g}^*.
\end{equation}
Imposing now $d\Omega\wedge\alpha$, we find
\[d(e^{136}+e^{145}+e^{235}-e^{246})\wedge e^7=0,\]
which by \eqref{eqn:hPlusV2} splits into
\[d(e^{136}+e^{145}-e^{246})\wedge e^7=0,\quad  de^{235}\wedge e^7=0\;;\]
from the second of which  \eqref{eqn:de2357} follows.
\end{proof}

Next we prove the ``if'' part of Theorem~\ref{thm:Classification}, by giving explicit examples for each solvable Lie group \textup{(\ref{eqn:h:abelian}\nobreakdash--\ref{eqn:h:L3:A0})}.

\emph{Abelian case}. Let  $a$, $b$ be real parameters with \[c=a^2 + b^2>0.\]
Every choice of $a$ and $b$ as above determines a semidirect product $\R^3\ltimes\R^4$, with an invariant basis of one-forms  satisfying
\begin{gather*}
de^1=0,\\
de^2=-\frac ac(a^{2}-2  b^{2})e^{56}+\frac bc(2 a^{2}-b^{2})(e^{15}+e^{26})+ a e^{34}-3 \frac{ b^{2} a }{c}e^{12}-\frac{1}{2}  {c} e^{17},\\
de^3= b e^{36}- b e^{45}+ a e^{13}+ a e^{24}+\frac{1}{2}  {c} e^{47},\\
de^4=0,\\
de^5=-\frac ac(a^{2}-2  b^{2})(e^{15}+e^{26})+ b e^{34}+\frac bc(2  a^{2}-b^{2})e^{12}-3 \frac{ b a^{2} }{c}e^{56}+\frac{1}{2}  {c} e^{67},\\
de^6=0,\\
de^7=-2 e^{12}-2 e^{34}-2 e^{56}.
\end{gather*}
Consider the  $\SU(3)$-structure determined by the choice of basis $e^1,\dotsc,e^7$. It is straightforward to verify that $d\alpha=-2F$ and $d\Omega\wedge\alpha=0$; it follows that this structure is contact and, by Theorem~\ref{thm:HypoContact}, the associated spinor is  generalized Killing. Applying Lemma~\ref{lemma:NablaOfForms}, we see that  in the chosen frame  the tensor $Q$ is given by the diagonal matrix
\begin{multline*}
\diag\biggl(1-\frac{1}{4} (a^{2}+ b^{2}),1+\frac{1}{4} (a^{2}+ b^{2}),1+\frac{1}{4} (a^{2}+ b^{2}),\\
1-\frac{1}{4} (a^{2}+ b^{2}),1+\frac{1}{4} (a^{2}+ b^{2}),1-\frac{1}{4} (a^{2}+ b^{2}),-3-\frac{3}{4} (a^{2}+ b^{2})\biggr).
\end{multline*}
\emph{Case \eqref{eqn:h:L4}}.\label{pg:eqn:h:L4}
This is a twofold case, since there is a choice of sign involved. Let $a$ be a real constant; consider the semidirect product $H\ltimes\R^4$ with structure constants determined by
\[ \left(0,0,- a e^{15},- a e^{16},a e^{13},a e^{14},-2 e^{12}-2 e^{34}-2 e^{56}\right).\]
By the same argument as before, the frame $e^1,\dotsc,e^7$ determines a contact $\SU(3)$-structure whose associated spinor is generalized Killing, and the tensor $Q$ is given by
\begin{gather*}
\diag(1,1,1,1,1,1,-3).
\end{gather*}
For the opposite sign, consider the group given by
\[ \left(0,-2  a^{2} e^{17}+2  a e^{36}-2  a e^{45},a e^{15},a e^{16},a e^{13}, a e^{14}, - 2 e^{12} - 2 e^{34} - 2 e^{56}\right).\]
The frame $e^1,\dotsc,e^7$ determines a contact $\SU(3)$-structure whose associated spinor is generalized Killing, and  $Q$ is determined by
\begin{gather*}
\diag(1-a^{2},1+a^{2},1,1,1,1,-3-a^{2}).
\end{gather*}

\emph{Case \eqref{eqn:h:L2}}.
Let $a$ be a real constant; consider the semidirect product  with structure constants determined by
\begin{multline*}
\biggl(0,-\frac{1}{2}  a^{2} e^{17}+ a e^{34}+ a e^{56},\frac{1}{2}  a^{2} e^{47}- a e^{24},a e^{14},\frac{1}{2}  a^{2} e^{67}- a e^{26},\\ a e^{16},-2 e^{12}-2 e^{34}-2 e^{56}\biggr).
\end{multline*}
The frame $e^1,\dotsc,e^7$ determines a contact $\SU(3)$-structure whose associated spinor is generalized Killing, and  $Q$ is determined by
\begin{gather*}
\diag\left(1-\frac{1}{4} a^{2},1+\frac{1}{4} a^{2},1+\frac{1}{4} a^{2},1-\frac{1}{4} a^{2},1+\frac{1}{4} a^{2},1-\frac{1}{4} a^{2},-3-\frac{3}{4} a^{2}\right).
\end{gather*}

\emph{Case \eqref{eqn:h:L3:A0}}.
Consider the semidirect product with structure constants determined by
\begin{multline*}
\biggl(
0, \frac{1}{2} e^{13}-\frac{3}{8} e^{17}+\frac{1}{2} e^{24}+e^{56}, \frac{1}{2} e^{12}-e^{34}+\frac{3}{8} e^{47}+\frac{1}{2} e^{56}, 0,\\ -e^{26}+\frac{1}{2} e^{36}-\frac{1}{2} e^{45}+\frac{3}{8} e^{67}, e^{16}, -2 e^{12}-2 e^{34}-2 e^{56}\biggr).
\end{multline*}
The frame $e^1,\dotsc,e^7$ determines a contact $\SU(3)$-structure whose associated spinor is generalized Killing, and  $Q$ is determined by
\begin{gather*}
\frac1{16}\diag\left(13,19,19,13,19,13,-57\right).
\end{gather*}

There are exactly two cases not appearing in the above examples, and so the proof of Theorem~\ref{thm:Classification} is reduced to the following non-existence result.
\begin{lemma}
Let $H$ be a solvable Lie group with structure equations  \eqref{eqn:h:nilpotent} or  \eqref{eqn:h:L3}. Then there is no representation of $H$ on $\R^4$ for which $H\ltimes V$ has a contact $\SU(3)$-structure whose associated spinor is generalized Killing.
\end{lemma}
\begin{proof}
Suppose $H\ltimes V$ has an $\SU(3)$-structure of the required type, and choose a basis of  invariant one-forms on $H\ltimes V$ as in Lemma~\ref{lemma:SemidirectStructure}. Since $\SO(3)\subset\SU(3)$ acts transitively on two-planes in $\R^3\cong\langle e^1,e^4,e^6\rangle$, we can assume that
\[\{\alpha\in \lie{h}^*\mid d\alpha \wedge\alpha=0\}\supset \Span{e^1,e^4}.\]
By acting further with an element of $\LieG{U}(1)\subset\SO(3)$, we can also assume that $de^1=0$. It follows that
\[de^1=0, \quad de^4=Ae^{14},\quad de^6=Be^{14}+Ce^{16},\]
where $A$ $B$, $C$ are constants. If $A$ is zero, then by hypothesis $H$ is nilpotent, and $C=0$; otherwise, $A$, $B$ and $C$ are all non-zero.

Let $\beta_j$ be a basis of the space \eqref{eqn:de2357}, and set $de^i=\sum_i a_{ij}\beta_j$. Imposing the linear conditions $dF=0$, $d\Omega\wedge\alpha=0$, we find
\begin{align*}
de^2&=  {(2 C-a_{39})} e^{56}+ a_{17} (e^{45}-e^{36})+ a_{13} (e^{15}+ e^{26})+ {(2 A-a_{25})}e^{34}\\
&- {(C+A-a_{39}-a_{25})} e^{12}+ a_{21} (e^{13}+ e^{24})+ a_{18} e^{47}+ a_{{1,12}} e^{67}+ a_{14} e^{17},\\
de^3&=  a_{{2,10}} (e^{45}-e^{36})- e^{24} a_{25}+ a_{21} e^{12}+ a_{28} e^{47}+ a_{{2,12}} e^{67}- a_{18} e^{17}\\
&- {(a_{21}-a_{{2,11}})} e^{34}- {(B+a_{17})} (e^{15}+ e^{26})+ {(A-a_{25})} e^{13}- a_{{2,11}} e^{56},\\
de^5&= - a_{17} e^{13}- a_{39} e^{26}- a_{{2,10}} e^{34}+ a_{13} e^{12}+ a_{{3,12}} e^{67}+ {(C-a_{39})} e^{15}+ a_{{2,12}} e^{47}\\
&- a_{{1,12}} e^{17}- {(B+a_{17})} e^{24}- {(a_{13}-a_{{2,10}})} e^{56} + a_{{2,11}} (e^{45}-e^{36}).
\end{align*}

Now, the condition $d^2=0$ determines an ideal $J$ of real polynomials in $a_{ij}$, $A$, $B$ and $C$. Lie groups with a structure of the required type correspond to points in the affine variety $V(J)$ determined by $J$. For the nilpotent case, we are interested in points with $A=C=0$. Calculations with \cocoa \cite{CocoaSystem} show that the polynomial $B^3$ lies in the ideal $J+(\{A,C\})$, and so $V(J)$ has no points with $A=C=0$, proving the statement in the nilpotent case.

For the case of \eqref{eqn:h:L3}, an analogous computation yields
\[a_{{2,12}} {\left((a_{14})^{2}+(a_{{2,12}})^{2}\right)}\in J.\]
It follows that $V(J+a_{2,12})=V(J)$. On the other hand, it turns out that
\[BC {(B^{2}+9 C^{2})} \in J+a_{2,12}.\]
This means that $V(J)$ has no points with both $B$ and $C$ different from zero, proving the rest of the statement.
\end{proof}

\section{Spin structures and submersions}
\label{sec:Submersions}A contact reduction is a two-step process, in which one first takes a submanifold, and then a quotient. The purpose of this section is to establish some formulae which will be needed to study the second step. It is in the context of Riemannian submersions that these formulae are presented most naturally.

To begin with, we study the relation among the Levi-Civita connections on the base and total space of a generic Riemannian submersion, using the language of principal bundles. Let $M^m$ be a manifold with an $\SO(m)$-structure $P_m$. Let $M^k$ be a manifold with an $\SO(k)$-structure
$P_k$, and let $\pi:M^m\to M^k$ be a Riemannian submersion.
 The tangent bundle of $M^m$ has an orthogonal splitting $\mathcal{H}\oplus\mathcal{V}$, where $\mathcal{V}=\ker\pi_*$. This defines a reduction of $P_m$ to  \[G=\SO(k)\times\SO(m-k)\;.\] Indeed, if we let $G$ act on $\R^m$ according to the splitting $\R^m=\R^k\oplus\R^{m-k}$,
the reduction $P_G$ is defined by
\[\mathcal{H}=P_G\times_G\R^k ,\quad \mathcal{V}=P_G\times_G\R^{m-k} . \]
We have a commutative diagram
\[\xymatrix{P_G\ar[r]^{d\pi}\ar[d] & P_k\ar[d]\\ M^m \ar[r]^\pi & M^k}\]
where $d\pi$ maps a frame $u:\R^m\to T_xM^m$ to a frame
\[d\pi(u)=\pi_*\circ u|_{\R^k}:\R^k\to T_{\pi(x)}M^k\;.\]
We say that a (local) section $s_m$ of $P_G$ is \dfn{$\pi$-related} to a section $s_k$ of $P_k$ if the diagram
\[\xymatrix{P_G\ar[r]^{d\pi} & P_k\\ M^m\ar[u]^{s_m} \ar[r]^\pi & M^k\ar[u]^{s_k}}\]
commutes. Then $s_ng$ is $\pi$-related to $s_k$ for all $\SO(m-k)$-valued functions $g$.

If $V$ is a $G$-module, an equivariant map $f_m:P_G\to V$ is \dfn{basic} if the diagram
\[\xymatrix{P_G\ar[rr]^{d\pi}\ar[dr]_{f_m} && P_k\ar[dl]^{f_k}\\&V}\]
commutes for some $f_k$. Then $f_k$ is   uniquely determined, and $\SO(k)$-equivariant; we say that $f_m$ is
\dfn{$\pi$-related} to $f_k$.

One can also  regard $f_m$  as a section $[s_m,x_m]$ of $P_G\times_G V$ where $x_m=f_m\circ s_m$. In this
language, a basic section $[s_m,\pi^*x_h]$ is $\pi$-related to $[s_h,x_h]$ if $s_m$ is $\pi$-related to
$s_h$. Moreover, all basic sections have this form. For instance, a  section of $\mathcal{H}$ is basic if it is  $\pi$-related to a vector field on $M^k$ in the usual sense.

\medskip
Let $\omega_k$ be the Levi-Civita connection form on $P_k$, so that the tautological form $\theta_k$ satisfies
\begin{equation}\label{eqn:dtheta}
d\theta_k+\omega_k\wedge\theta_k=0\;.
\end{equation}
The tautological form on $P_m$ restricted to $P_G$ can be written as $\theta_m=(\theta_h,\theta_v)$. By
construction $\theta_h=(d\pi)^*\theta_k$. Then \eqref{eqn:dtheta} yields
\begin{equation}
\label{eqn:dthetah}
0=(d\pi)^*(d\theta_k+\omega_k\wedge\theta_k)=d\theta_h+(d\pi)^*\omega_k\wedge\theta_h\;.
\end{equation}
On the other hand, the connection form $\omega_m$ restricted to $P_G$ is an $\so(m)$-valued $1$-form, which we
can decompose into blocks as
\[\omega_m=\begin{pmatrix}\omega_h & -A^T \\ A & \omega_v\end{pmatrix}.\]

\begin{lemma}
\label{lemma:ONeill}
If $X$ is a ``horizontal'' vector field on $P_G$, in the sense that $\theta_v(X)$ is zero, then \[\omega_h(X)=(d\pi)^*\omega_k(X)\;.\]
\end{lemma}
\begin{proof}
At points of $P_G$, the analogue of \eqref{eqn:dtheta} yields
\[0= (d\theta_h + \omega_h\wedge\theta_h -A^T\wedge \theta_v , d\theta_v +A\wedge\theta_h+\omega_v\wedge\theta_v).\]
Comparing with \eqref{eqn:dthetah}, it follows that
\begin{equation}
\label{eqn:hVSk} (\omega_h-(d\pi)^*\omega_k)\wedge\theta_h=A^T\wedge\theta_v.
\end{equation}
This means that $\omega_h$ is not $\pi$-related to $\omega_k$. However, consider the standard isomorphism
\[\partial:\R^k\otimes\Lambda^2\R^k\to\Lambda^2\R^k\otimes\R^k.\]
Equation \eqref{eqn:hVSk} tells us that, if we define  a map
\[\beta:P_G\to \R^k\otimes\Lambda^2\R^k,\quad \langle\beta,\theta_h(X)\rangle=(\omega_h-(d\pi)^*\omega_k)(X),\]
then $\partial \beta=0$ and so $\beta=0$. Thus, if $X$ is horizontal then
\[(\omega_h-(d\pi)^*\omega_k)(X)=0\;.\qedhere\]
\end{proof}
\begin{remark}
To illustrate the meaning of $A$, apply the formula \[\nabla^m_X [s,x]=\left[s,s^*\omega_m(X) x + dx(X)\right]\;;\]
for a horizontal vector field $Y$, which we write as a section $[s,h]$ of $\mathcal{H}$, we find that the component in $\mathcal{V}$ of its covariant derivative satisfies
\[(\nabla^m_X Y)_v=[s,s^*A(X) h]\;.\]
Thus, $A$ is the principal bundle version of the  O'Neill tensor.
\end{remark}

\medskip
In order to pass from $\SO(n)$-structures to spin structures, consider the chain of inclusions
\[\SO(k)\to \SO(k)\times\SO(m-k)\to \SO(m)\;,\]
which gives rise to a commutative diagram
\begin{equation}\label{eqn:Diagram}\begin{gathered}
\xymatrix{ \Cl(m) & \ar[l] \Spin(m)\ar[r]^{\Ad} & \SO(m)\\
\Cl(k)\ar[u]^{j_k^m}&  \Spin(k)\ar[l]\ar[u]^{j_k^m}\ar[r]^{\Ad} & \SO(k)\ar[u]^{j_k^m}
}
\end{gathered}
\end{equation}
Explicitly, if $e_1,\dots,e_m$ is the standard basis of $\R^m$, we define the homomorphism of Clifford algebras $j_k^m$ by
\[ j_k^m(e_i)=e_{k+1}\cdot e_i\;.\]
The homomorphism is not unique, but this choice  has the advantage that it makes the two half-spin representations $\Sigma^{2n}_\pm$ into (irreducible) $\Cl(2n-1)$-modules. A more relevant consequence of the definition is the equality
 \begin{equation*}
 j_l^m\circ j_k^{l}=j_k^m \text{ for } k<l<m\;.
\end{equation*}

\smallskip
We can now return to Riemannian submersions, and introduce spinors in the picture. The first problem is defining the spin structure on $M^k$ in terms of the spin structure on $M^m$. We shall assume from now on that $\mathcal{V}$ is parallelizable. Then we can replace $G$ with $\SO(k)$ in the
construction; indeed, in the rest of this section we set
\[G=\SO(k)\subset\SO(m)\;.\]
Suppose that $M^m$ is spin, and fix a spin structure on $M^m$, i.e. a principal bundle $P_{\Spin}$ with structure group $\Spin(m)$ and a $2:1$ equivariant projection $\Ad:P_{\Spin}\to P_{m}$.
We define a principal bundle with fibre $\Spin(k)$
\[P=\Ad^{-1}(P_G)\]
on $M^m$. The diagram \eqref{eqn:Diagram} shows that $P$ is a good candidate for the pullback to $M^m$ of a spin structure on $M^k$. We shall now make this notion precise.

Denoting the generic point of $P_G$ by $(x;u)$, where $x$ is in $M^m$ and $u$ is a point in the fibre of $x$, the manifold $P_k$ can be identified with the quotient $P_G/\sim$, where $(x_0;u_0)\sim (x_1;u_1)$ if and only if $\pi(x_0)=\pi(x_1)$ and there is a basic section $s$ such that $s(x_i)=u_i$. By construction, sections of $P_G/\sim$ can be identified with basic sections of $P_G$. Similarly, we say that a section $\tilde s$ of $P$ is basic if $\Ad(\tilde s)$ is basic, and define an equivalence relation on $P$ by $(x_0;\gamma_0)\sim (x_1;\gamma_1)$ if and only if $\pi(x_0)=\pi(x_1)$ and there is a basic section $\tilde s$ such that $\tilde s(x_i)=\gamma_i$.
\begin{lemma}
\label{lemma:SpinorBundleOnMk}
$P/\sim$ is a spin structure on $M^k$. Accordingly, a spinor on $M^k$ is given by a basic section of $P\times_{\Spin(k)}\Sigma^k$, and its Levi-Civita covariant derivative is given by
\[\nabla [\tilde s,\pi^*\psi]=[\tilde s,\tilde s^*(d\pi)^*(\omega_k \psi) +\pi^*d\psi]\;,\]
where $\tilde s$ is a local basic section of $P$ and $\psi\colon M^k\to\Sigma^k$.
\end{lemma}
\begin{proof}
For the first part, it suffices to show that $\Ad:P/\sim\to P_G/\sim$ is a two-sheeted covering. Suppose that
$\Ad(x_0;\gamma_0)\sim \Ad(x_1;\gamma_1)$. Then there exists a section $s$ of $P_G$ with $s(x_i)=\Ad(\gamma_i)$. Since $P_G$ is trivial on the fibres of $\pi$, we can lift $s$ to a unique section $\tilde s$ of $P$ with $\tilde s(x_0)=\gamma_0$. Hence $\Ad(\tilde s(x_i))=\Ad(\gamma_i)$, and
\[(x_0;\gamma_0)\sim (x_1;\gamma_1)\iff \tilde s(x_i)=\gamma_i\;.\]
The pullback  $(d\pi)^*\omega_k$ is clearly a torsion-free, $\so(k)$-valued connection form, and therefore coincides with the Levi-Civita connection.
\end{proof}
\begin{remark}
The operator $\nabla_X$ defined in Lemma~\ref{lemma:SpinorBundleOnMk} is defined for arbitrary $X$ on $M^m$. It represents the Levi-Civita covariant derivative on $M^k$ when $X$ is the basic lift of a vector field on $M^k$; when $X$ is vertical, $\nabla_X$ is zero.
\end{remark}

The spinor bundle on $M^m$ can be identified with
\[P\times_{\Spin(k)} \Sigma^m\;,\]
where the action of $\Spin(k)$ on $\Sigma^m$ is induced by the map $j_k^m$ of Diagram~\ref{eqn:Diagram};
the connection form $\omega_m$ defines a covariant derivative operator $\nabla^m$ on this vector bundle. The operator $\nabla$ of Lemma~\ref{lemma:SpinorBundleOnMk} can clearly be defined for basic sections of any bundle associated to $P$; in particular, we can compare $\nabla$ and $\nabla^m$ on basic sections of $P\times_{\Spin(k)} \Sigma^m$. This will be the main result of this section, which will only be stated for $k=m-1$, since we are ultimately interested in $\LieG{U}(1)$ reductions. This assumption, by which nothing would have been gained up to this point, has two useful consequences: we can represent $A(X)$ by a $\mathcal{H}$-valued one-form, namely
\[\underline{A}=[s,s^*A_{j}e_j],\]
and $\omega_v$ vanishes. To state the result, we need to introduce the vector bundle map
\[\underline{j}\colon P\times_{\Spin(k)} \Cl(k)\to P\times_{\Spin(k)} \Cl(m),\]
induced by $j_k^m=j_{m-1}^m$.
\begin{proposition}
\label{prop:Submersion}
Let $\tilde s$ be a basic section of $P$ and $\psi_m\colon M^k\to\Sigma^m$. Then
\begin{equation*}
\nabla^m_X [\tilde s,\pi^*\psi_m]=\nabla_X [\tilde s,\pi^*\psi_m] -\frac12 \underline{j}(\underline{A})(X)\cdot [\tilde s, \pi^*\psi_m]\;.
\end{equation*}
\end{proposition}
\begin{proof}
 By Lemma~\ref{lemma:ONeill} and Lemma~\ref{lemma:SpinorBundleOnMk}, if $X$ is
horizontal
\begin{multline*}
\nabla_X [\tilde s,\psi_m] = [\tilde s, \tilde s^* (d\pi)^*(\omega_k\psi_m)(X) + \pi^*d\psi_m(X)]=\\
=[\tilde s, \tilde s^*\omega_h(X)\pi^*\psi_m + \pi^*d\psi_m(X)]\,.
\end{multline*}
On the other hand, recall the expression of $\omega_m$, which can be rewritten as
\[\omega_m=\omega_h + \sum_{1\leq j< m} A_{j}(e_{mj}-e_{jm})\;;\]
the second component acts on spinors as
\[-\frac12 \sum_{j} A_{j}e_m\cdot e_j\;.\]
Thus, for the covariant derivative $\nabla^m$ we obtain
\begin{multline*}\nabla^m_X\left[\tilde s,\pi^*\psi_m\right]=\left[\tilde s, \tilde s^*\omega_m(X)\pi^*\psi_m + \pi^*d\psi_m(X)\right]\\
 = \left[\tilde s,\tilde s^*\omega_h(X)\pi^*\psi_m  -\frac12 s^*A_{j}(X) e_m\cdot e_j\cdot\psi_m +
\pi^*d\psi_m(X)\right].
\end{multline*}
The statement now follows from the definition of $\underline{j}$.
\end{proof}

\section{Contact reduction and spinors}
\label{sec:Contact}
This section contains the main result of this paper. We return to the situation of Theorem~\ref{thm:HypoContact}, with a slight change of language: now, we  regard the contact structure and metric as fixed, and identify  a rank one complex bundle of ``compatible'' spinors, determined by the $\LieG{U}(n)$-structure. The reduction from  $\LieG{U}(n)$ to $\SU(n)$ is represented by the choice of a unit section of this bundle. We show that, in the presence of a circle action that preserves the contact metric structure and a compatible spinor, one can define a compatible spinor on the contact reduction. Applying Proposition~\ref{prop:Submersion}, we are able to determine sufficient conditions for the generalized Killing condition to be preserved by the reduction process. In the language of Theorem~\ref{thm:HypoContact}, at least in the real analytic category, this means that the symplectic reduction of a symplectic cone with a Calabi-Yau metric is again Calabi-Yau, under  certain conditions.

Let $M$ be a $2n+1$-dimensional manifold. A $\LieG{U}(n)$-structure $P_{\LieG{U}}$ on $M$ identifies differential forms $\alpha$, $F$ by \[\alpha=e^{2n+1},\quad F=e^{12}+\dots+e^{2n-1,2n}\;;\]
like in Section~\ref{sec:GKSpinors}, we say that the structure is \dfn{contact} if $d\alpha=-2F$.  In this case, we shall also refer to $P_{\LieG{U}}$ as a \dfn{contact metric structure}.

Now suppose that $M$ is spin; let $P_{\Spin}$ be a spin structure compatible with the given metric and orientation. The preimage of $\LieG{U}(n)$ under $\Ad$ is a connected subgroup  $\tilde{\LieG{U}}(n)\subset\Spin(2n+1)$. Passing to the principal bundles,  the preimage of the $\LieG{U}(n)$-structure under $\Ad$ is a  $\tilde{\LieG{U}}(n)$-reduction of $P_{\Spin}$; we shall denote it by $P_{\tilde{\LieG{U}}}$. Thus, a spinor on $M$ is a section of
\[P_{\tilde{\LieG{U}}}\times_{\tilde{\LieG{U}}(n)}\Sigma^{2n+1}\;.\]
By \cite{FriedrichKim:TheEinsteinDiracEquation}, $\Sigma^{2n+1}$ splits as $\Sigma_0^{2n+1}\oplus\dots\oplus\Sigma_n^{2n+1}$ as a representation of $\tilde{\LieG{U}}(n)$, where in particular
\begin{equation}
\label{eqn:Sigma0}
\psi\in\Sigma_0^{2n+1}\;\iff\; e^{2n+1}\cdot \psi =i^{2n+1}\psi\;,\quad e^{2k-1,2k}\cdot\psi =-i\psi.
\end{equation}

Thus, we have identified the bundle of \dfn{compatible} spinors
\begin{equation}
\label{eqn:SpinorsOfTypeSigma0}
P_{\tilde{\LieG{U}}}\times_{\tilde{\LieG{U}}(n)}\Sigma_0^{2n+1},
\end{equation}
whose sections $\psi$ are
characterized by
\[\alpha\cdot \psi =i^{2n+1}\psi\;,\quad F\cdot\psi =-ni\psi\;.\]

Assume now that $S^1$ acts  on $M$ preserving both metric and contact form, so  that the fundamental vector field $X$
satisfies
\[\Lie_X\alpha=0=\Lie_X F\;.\]
The moment map is defined by
\[\mu=\alpha(X)\;;\]
we assume that $0$ is a regular value of $\mu$. Denote by $\iota:M_0\to M$ the hypersurface $\mu^{-1}(0)$. The contact reduction of $M$ is by definition \cite{Geiges:Constructions}
\[M//S^1 = M_0/S^1\;.\]
By construction, $\alpha(X)=0$ on $M_0$, and the
tangent bundle of $M_0$ consists of vectors $Y$ with
\begin{equation}
\label{eqn:TM0}
0=d(\alpha(X))Y=-(X\hook d\alpha)Y=-2d\alpha(X,Y)=F(X,Y) \;.
\end{equation}
In the pullback bundle $\iota^*P_{\LieG{U}}$, consider the set $P_{\LieG{U}(n-1)}$ of those frames $u$
such that
\[[u,e_{2n-1}]=tX_{p(u)}\;,\]
where $t$ is a positive function. This defines a $\LieG{U}(n-1)$-structure on $M_0$, inducing in in turn a contact $\LieG{U}(n-1)$-structure on $M//S^1$. Moreover, we can define a unit normal vector field $\nu$, by
\[\nu_{p(u)}=[u,e_{2n}], \quad u\in P_{\LieG{U}(n-1)}\;.\]
We can also think of $\nu$ as the vector field dual to the one-form $t^{-1}X\hook F$, and write
\[\nu=(t^{-1}X\hook F)^\sharp.\]

We now show that the choice of an invariant compatible spinor on $M$ determines a compatible spinor on $M//S^1$.
By \eqref{eqn:TM0}, the inclusion of $\R^{2n-1}$ into $\R^{2n+1}$ that determines the $\LieG{U}(n-1)$-structure on $M//S^1$ corresponds to the basis
\begin{equation}
 \label{eqn:subbasis}
e_1,\dotsc,e_{2n-2},e_{2n+1}\;.
\end{equation}
Accordingly, the algebra homomorphism $j_{2n-1}^{2n+1}$ of \eqref{eqn:Diagram} is given by
\[j(e_k)=e_{2n-1}\cdot e_k\;.\]
Using $j$, we can view $\Sigma^{2n+1}$ as a $\Cl(2n-1)$-module that splits into two irreducible components, one of which is isomorphic to $\Sigma^{2n-1}$; it can be identified by the action of the volume element.
\begin{lemma}\label{lemma:Sigma0}
The space of spinors in $\Sigma^{2n-1}\subset\Sigma^{2n+1}$ that are compatible with the  \mbox{$\LieG{U}(n-1)$}\nobreakdash-\hspace{0pt}structure determined by the frame \eqref{eqn:subbasis}
is given by
\begin{align*}
\Sigma_0^{2n-1}=\{\psi-e_{2n-1}\cdot\psi,\;\psi\in\Sigma_0^{2n+1}\}.
\end{align*}
\end{lemma}
\begin{proof}
If a spinor in $\Sigma^{2n+1}$ satisfies the compatibility conditions \eqref{eqn:Sigma0} with respect to to the \mbox{$\LieG{U}(n-1)$}\nobreakdash-\hspace{0pt}structure, then in particular it belongs to \mbox{$\Sigma^{2n-1}\subset\Sigma^{2n+1}$}.
Therefore, it suffices to check that for $\psi$ in $\Sigma_0^{2n+1}$
\[j(e_{2n+1})\cdot(\psi-e_{2n-1}\cdot\psi) = e_{2n-1}\cdot e_{2n+1}\cdot (\psi-e_{2n-1}\cdot\psi)=i^{2n-1}(\psi-e_{2n-1}\cdot\psi),\]
and, for $1\leq k<n$,
\[
j(e^{2k-1,2k})\cdot(\psi-e_{2n-1}\cdot\psi)=(1-e_{2n-1})\cdot e^{2k-1,2k}\cdot\psi=
-i(\psi-e_{2n-1}\cdot\psi).\qedhere
\]
\end{proof}
Now define a principal bundle $P_{\LieG{\tilde U}(n-1)}$ on $M_0$ so that
\[
\xymatrix{ P_{\LieG{\tilde U}(n-1)} \ar[r]\ar[d] & \iota^*P_{\LieG{\tilde U}}\ar[d]\\
P_{\LieG{U}(n-1)} \ar[r] & \iota^*P_{\LieG{U}}}
\]
is a commuting diagram of principal bundles on $M_0$, and all maps are equivariant.
A section $\psi$ of \eqref{eqn:SpinorsOfTypeSigma0} pulls back to a section $\iota^*\psi$ of
\[P_{\tilde{\LieG{U}}(n-1)}\times_{\tilde{\LieG{U}}(n-1)}\Sigma^{2n+1}_0\;.\]
Now assume that $\psi$ is $S^1$-invariant. With Lemma~\ref{lemma:Sigma0} in mind, we define \mbox{$\psi^\pi = \iota^*\psi - t^{-1}X\cdot\iota^*\psi$}, which we
rewrite as
\begin{equation*}
\psi^\pi=\iota^*\psi + i\nu\cdot\iota^*\psi\;.
\end{equation*}
By Section~\ref{sec:Submersions}, the spin structure and the spinor $\psi^\pi$ on $M_0$ induce a spin
structure and  spinor on $M_0/S^1$. By Lemma~\ref{lemma:Sigma0}, this spinor is compatible with the contact
metric structure of $M_0/S^1$. We can now state our main result.
\begin{theorem}
\label{thm:HypoToHypo}
Let $M$ be a manifold of dimension $2n+1$ with a contact $\LieG{U}(n)$-structure $(g,\alpha, F)$ and a
compatible generalized Killing spinor $\psi$, so that
\[\nabla_Y^{2n+1}\psi=\frac12 Q(Y)\cdot\psi\;, \quad Y\in TM\;,\]
where $Q$ is a symmetric endomorphism of $TM$. Suppose furthermore that $S^1$ acts on $M$ preserving both  structure and spinor, zero is a regular value for the moment map $\mu$, and $S^1$ acts freely on $\mu^{-1}(0)$. Then the  spinor $\psi^\pi$ induced on $M//S^1$ is compatible with the induced contact metric structure; in addition, $\psi^\pi$  is generalized Killing  if and only if at each point of $\mu^{-1}(0)$
\[ dt\in \Span{X\hook F,\alpha},\]
where $X$ is the fundamental vector field associated to the $S^1$ action, and  $t$ its norm. In this case,
\[\nabla_Y^{2n-1} \psi^\pi = \frac12 \underline{j}(B(Y))\cdot\psi^\pi\;, \quad Y\in T(M//S^1),\]
where $B$ is the symmetric endomorphism of $T(M//S^1)$ given by
\[B(Y,Z)=-Q(Y,Z) -Q(t^{-1}X,t^{-1}X)\alpha(Y)\alpha(Z).\]
\end{theorem}
In the statement of Theorem~\ref{thm:HypoToHypo}, $\nabla^{2n+1}$ denotes the covariant derivative on $M$, and $\nabla^{2n-1}$ the covariant derivative on $M//S^1$ defined in Lemma~\ref{lemma:SpinorBundleOnMk}; we shall also consider the covariant derivative $\nabla$ on $M_0$.
Our calculations will also involve the Weingarten tensor $W$ of $M_0\subset M$, and  the $\mathcal{H}$-valued one-form $A$  determined by the Levi-Civita connection on $M_0$ (see Proposition~\ref{prop:Submersion}).
Finally, from now on $Y$ represents a generic basic vector field on $M_0$. We can now establish some useful formulae.
\begin{lemma}
\label{lemma:LemmaHypoToHypo}
In the hypotheses of Theorem~\ref{thm:HypoToHypo} (with no assumptions on $\psi^\pi$), at each point of $M_0$
\begin{gather}
\label{eqn:Rnu}
(-1)^n Q(Y,\nu)\alpha+ A(Y) +W(Y)\hook F\in\R\nu\;,\\
\label{eqn:Automatic}Q(Y,t^{-1}X)=(-1)^n\alpha(W(Y)),\\
\label{eqn:WX} W(X)=2(dt)^\sharp\hook F + (-1)^n  Q(X,t^{-1}X)\alpha.
\end{gather}
\end{lemma}
\begin{proof}
By definition, the tensor $A$ is characterized by
\[\langle A(Y),Z\rangle = \langle \nabla_Y Z,t^{-1}X\rangle =2F(\nabla_Y^{2n+1} Z,\nu)\]
for all horizontal vector fields $Y,Z$ on $M_0$. By Lemma~\ref{lemma:NablaOfForms}, it follows that
\begin{multline*}
0=\Lie_Y F(Z,\nu)=(\nabla^{2n+1}_Y F)(Z,\nu)+F(\nabla^{2n+1}_Y Z,\nu)+F(Z,\nabla^{2n+1}_Y\nu)=\\
=\frac12\bigl((-1)^n\alpha(Z)Q(Y,\nu)+\langle A(Y),Z\rangle +2F(Z,-W(Y))\bigr).
\end{multline*}
This expression also vanishes trivially for $Z=X$, and so  \eqref{eqn:Rnu} follows.

To prove the second equation, one uses Lemma~\ref{lemma:NablaOfForms} and the fact that $Q$ is symmetric, from which
\[W(Y,\alpha^\sharp)=\langle \nabla^{2n+1}_Y\alpha,\nu\rangle =(-1)^n 2F(Q(Y),\nu)\;,\]
which is equivalent to  \eqref{eqn:Automatic}.

Finally, one can write
\[W(X)=-\nabla^{2n+1}_{X} (t^{-1}X\hook F)=\left(-\nabla^{2n+1}_{X} t^{-1}X\right)\hook F + (-1)^n  Q(X,t^{-1}X)\alpha;\]
on the other hand, one has
\begin{multline}
0=\Lie_X\langle t^{-1}X,\cdot\rangle  - \langle \nabla^{2n+1}_X (t^{-1}X),\cdot\rangle -\langle t^{-1}X,\nabla^{2n+1}_X \rangle\\
= - \langle \nabla^{2n+1}_X t^{-1}X,\cdot\rangle-\langle t^{-1}X,\nabla^{2n+1} X\rangle,
\end{multline}
whence
$\nabla^{2n+1}_X (t^{-1}X) = -2(dt)^\sharp$, giving \eqref{eqn:WX}.
\end{proof}
We can now prove Theorem~\ref{thm:HypoToHypo}.
\begin{proof}[Proof of Theorem~\ref{thm:HypoToHypo}]
By the analogue of Proposition~\ref{prop:Submersion} for immersions (see e.g. \cite{BarGauduchonMoroianu}), the covariant derivatives on $M_0$ and $M$ are related by
\[\nabla_Y \iota^*\psi = \iota^*(\nabla_Y^{2n+1}\psi)+\frac12 \nu\cdot W(Y)\cdot\iota^*\psi\;.\]
Hence, applying this formula to $\psi^\pi$,
\begin{multline*}
\nabla_Y\psi^\pi=\iota^*\left(\nabla_Y^{2n+1}(\psi+i\nu\cdot\psi)\right)+\frac12 \nu\cdot W(Y)\cdot
\iota^*(\psi+i\nu\cdot\psi)\\
=(1+i\nu)\cdot\left(\nabla_Y^{2n+1}\psi-\frac i2 W(Y)\cdot\psi\right),
\end{multline*}
where we have used $\nabla_Y^{2n+1}\nu=-W(Y)$. By Proposition~\ref{prop:Submersion}, denoting by $A$ the $\mathcal{H}$-valued one-form on $M^0$ determined by the Levi-Civita connection on $M_0$,
\begin{multline*}
\nabla^{2n-1}_Y\psi^\pi=\nabla_Y\psi^\pi+\frac12 t^{-1}X\cdot A(Y)\,\cdot\psi^\pi
=(1+i\nu)\cdot\left(\nabla_Y^{2n+1}\psi-\frac 12 (A(Y)+iW(Y))\cdot\psi\right),
\end{multline*}
and so by hypothesis
\begin{equation}
\label{eqn:OnReduction} \nabla^{2n-1}_Y\psi^\pi=\frac12(1+i\nu)\cdot(Q(Y)-A(Y)-iW(Y))\cdot\psi.
\end{equation}
The induced spinor $\psi^\pi$ is generalized Killing  on $M//S^1$ if and only if
\[\nabla_Y^{2n-1}\psi^\pi=\frac12\underline{j}(B(Y))\cdot\psi^\pi\]
for $B$ symmetric, or equivalently,
\[\nabla_Y^{2n-1}\psi^\pi=\frac12t^{-1}X\cdot B(Y)\cdot(1+i\nu)\cdot\psi=-\frac12(1+i\nu)\cdot B(Y)\cdot\psi\;.\]
Comparing with \eqref{eqn:OnReduction}, we find that $\psi^\pi$ is generalized Killing if and only if
\begin{equation}
\label{eqn:B}
(1+i\nu)\cdot(B(Y)+Q(Y)-A(Y)-iW(Y))\cdot\psi=0 \text{ for some symmetric }B\;.
\end{equation}
Now
\[
(W(Y)\hook F)\cdot\psi=iW(Y)\cdot\psi+(-1)^n\alpha(W(Y))\psi\;,
\]
so $B(Y)$ is characterized by
\[ (1+i\nu)(B(Y)+Q(Y)-A(Y)-W(Y)\hook F+(-1)^n\alpha(W(Y)))\cdot\psi=0,
\]
which we can split into components as
\begin{multline*}
 (B(Y)+Q(Y)-A(Y)-W(Y)\hook F)\cdot\psi\\
-(B(Y)+Q(Y)-A(Y)-W(Y)\hook F)\cdot\nu\cdot i\psi\\
+(-1)^n\alpha(W(Y))\psi\\
+(-2Q(Y,\nu)+2W(Y,t^{-1}X))i\psi\\
+(-1)^n\alpha(W(Y))\nu\cdot i\psi =0.
\end{multline*}
This is only possible if all four summands are in $\SpanC{\nu\cdot\psi,\psi}$, and equivalent to requiring that at each point (using also \eqref{eqn:Automatic})
\begin{multline}
\label{eqn:BOAW}
B(Y)+Q(Y)-A(Y)-W(Y)\hook F \\=Q(Y,t^{-1}X) t^{-1}X   +
(Q(Y,\nu)-W(Y,t^{-1}X))((-1)^n\alpha + \nu).
\end{multline}
By \eqref{eqn:Rnu}, we obtain
\begin{multline*}
B(Y)+Q(Y) +(-1)^nQ(Y,\nu)\alpha - Q(Y,t^{-1}X) t^{-1}X \\
-(Q(Y,\nu)-W(Y,t^{-1}X))((-1)^n\alpha + \nu)\in\R\nu,
\end{multline*}
which is equivalent to
\begin{gather}
\label{eqn:DefineB} B(Y)=-Q(Y)^\perp-(-1)^nW(Y,t^{-1}X)\alpha\;,
\end{gather}
where the superscript `$\perp$' denotes projection on the orthogonal complement of $\Span{X,\nu}$.
In light of \eqref{eqn:DefineB},  we can rewrite \eqref{eqn:BOAW} as
\[
-A(Y)-W(Y)\hook F -    (-1)^n Q(Y,\nu)\alpha+W(Y,t^{-1}X)\nu=0,
\]
which by Lemma~\ref{lemma:LemmaHypoToHypo} is automatically satisfied.

Summing up, \eqref{eqn:B} is equivalent to the tensor $B$ defined in \eqref{eqn:DefineB} being symmetric.
This is equivalent to $W(Y,t^{-1}X)$ being zero whenever $Y$ is horizontal and $\alpha(Y)$ is zero, or, by the symmetry of $W$, to
\[W(X)\in\Span{\alpha^\sharp,X}.\]
By Lemma~\ref{lemma:LemmaHypoToHypo}, the first part of the statement follows. The second part is now a straightforward consequence of  \eqref{eqn:DefineB} and Lemma~\ref{lemma:LemmaHypoToHypo}.
\end{proof}

\begin{remark}
The statement of Theorem~\ref{thm:HypoToHypo} is essentially local. In fact, if one replaces the Lie group $S^1$ with $\R$ the proof carries through, provided the contact reduction $M//\R$ is well defined and smooth.
\end{remark}

In our language, an $\alpha$-Einstein-Sasaki structure on $M^{2n+1}$, $n>1$ can be characterized as a contact $\LieG{U}(n)$-structure admitting a compatible generalized Killing spinor with
\[Q(Y)=a\, Y +b\,\alpha(Y)\alpha^\sharp,\]
where $a$ and $b$ are constants (see \cite{FriedrichKim:TheEinsteinDiracEquation}). By \cite{GrantcharovOrnea}, the contact reduction of an $\alpha$\nobreakdash-Einstein-Sasaki structure is Sasaki. As a consequence of Theorem~\ref{thm:HypoToHypo}, we obtain the following:
\begin{corollary}
 \label{cor:alphaEinsteinSasaki}
Let $M$ be a manifold of dimension $2n+1$ with an $\alpha$\nobreakdash-Einstein-Sasaki structure $(g,\alpha, F,\psi)$, and let $S^1$ act on $M$ preserving the structure in such a way that $0$ is a regular value for the moment map $\mu$ ans $S^1$ acts freely on $\mu^{-1}(0)$. Then the Sasaki quotient $M//S^1$ is also $\alpha$-Einstein if and only if
\[dt\in \Span{ X\hook F,\alpha}\]
 at each point of $\mu^{-1}(0)$,  where $X$ is the fundamental vector field associated to the $S^1$ action and $t$ is its norm.
\end{corollary}

\section{Examples}
\label{sec:Examples}
In this section we apply Theorem~\ref{thm:HypoToHypo} to two concrete examples in dimension seven, obtaining hypo-contact structures in dimension five; one of the resulting structures is the nilpotent example appearing in \cite{DeAndres:HypoContact}, and the other is new.
\subsection{The Heisenberg group}
As mentioned in the introduction, the Heisenberg group $G_{2n+1}$ of dimension $2n+1$ has an $\alpha$-Einstein-Sasaki structure (see also \cite{TomassiniVezzoni}). We can represent $G_{2n+1}$ by a basis of left-invariant one-forms $e^1,\dotsc,e^{2n+1}$ satisfying
\[de^1=0,\dotsc, de^{2n}=0, de^{2n+1}=-2\left(e^{12}+\dotsb+e^{2n-1,2n}\right).\]
The choice of a basis $e^1,\dotsc, e^{2n+1}$ determines an $\SU(n)$-structure by \eqref{eqn:StandardFrame}, which is $\alpha$-Einstein-Sasaki with
\[Q(e_i,e_i)=(-1)^{n+1}, Q(e_{2n+1},e_{2n+1})=(-1)^nn,\]
and the other components of $Q$ equal to  zero.
Now let $X$ be the right-invariant vector field with $X_e=e_{2n-1}$. The Lie group $\{\exp ue_{2n-1}\}$ is closed in $G$, and acts on $G$ on the left, preserving the $\SU(n)$-structure, with associated fundamental vector field $X$. By construction,
\[\mu(g)=0\iff L_{g*}^{-1}X_g\in\langle e_{2n-1}\rangle,\]
so $\mu^{-1}(0)$ is the subgroup with Lie algebra $\langle e_1,\dotsc,e_{2n-1},e_{2n+1}\rangle$, which contains $\{\exp ue_{2n-1}\}$ as a normal subgroup. It follows that the contact reduction $G_{2n+1}//S^1$ is the Heisenberg group $G_{2n-1}$, consistently with Theorem~\ref{thm:HypoToHypo}.
\begin{remark}
 This is not the only way one can obtain the Heisenberg group by means of a reduction;  for example, the above construction applies equally well to  the examples of Section~\ref{sec:ExamplesNoQuotient} relative to the Lie group $H$ with structure constants \eqref{eqn:h:L4}, setting $X_e=e_2$. On the other hand, quotients of  semidirect products $H\ltimes V$ do not generally satisfy Theorem~\ref{thm:HypoToHypo} when $H$ is solvable.
\end{remark}

\subsection{A new compact example}
Consider the semidirect product $\SU(2)\ltimes_\phi\h$, where $\phi\colon\SU(2)\to \GL(\h)$ is given by quaternionic multiplication. By \cite{FinoTomassini}, this Lie group has a left-invariant weakly integrable generalized $\Gtwo$-structure; in fact, we shall see that this structure satisfies Theorem~\ref{thm:HypoToHypo}, giving rise to a hypo-contact structure on $S^2\times T^3$.

We shall use quaternionic coordinates, following \cite{Salamon:Tour}, and write the generic element of \mbox{$\SU(2)\ltimes\h$} as $(p,q)$. There are two natural left-invariant quaternionic forms, which give rise to a left-invariant basis of one-forms $e^1,\dotsc,e^7$ by
\[p^{-1}dp=-2ie^1-2je^4 -2ke^6\;,\quad p^{-1}dq=e^3+ie^5+je^7+ke^2\;.\]
Notice that $p^{-1}=\overline{p}$. Also, conjugation and right quaternion multiplication are equivariant, and
\[(d\overline p)p = \overline{p^{-1}dp}, \quad p^{-1}dq\,j= -e^7-ie^2+je^3+ke^5\]
are also invariant forms.
From the identities
\begin{align*}
d(p^{-1}dp) &=d\overline{p}\wedge dp=(d\overline{p})p\wedge p^{-1} dp,\\
d(p^{-1}dq) &=d\overline{p}\wedge dq=(d\overline{p})p\wedge p^{-1}  dq,
\end{align*}
we deduce that the structure constants are given by
\begin{multline*}
\label{eqn:AnnaAdriano}
\left((4e^{46},2(-e^{36}-e^{45}+e^{17}),2(-e^{15}+e^{26}-e^{47}),-4e^{16},\right.\\
\left.2(e^{13}-e^{24}-e^{67}),4e^{14},-2(e^{12}+e^{34}+e^{56})\right).
\end{multline*}
Again, the basis  $e^1,\dotsc,e^7$ defines a contact $\SU(3)$-structure by \eqref{eqn:StandardFrame}, and an identification of left-invariant spinors with $\Sigma\cong\C^8$; the associated spinor $\psi=u_0$
satisfies \eqref{eqn:GKSpinor} with
\[Q=\diag(2,0,0,2,0,2,0).\]
In quaternionic terms, denoting by $\odot$ the symmetric product, we can describe the metric tensor as
\[g=\frac14\overline{p^{-1}dp}\odot p^{-1}dp + \overline{p^{-1}dq} \odot p^{-1}dq=\frac14d\overline{p}\odot dp+d\overline{q}\odot dq,\]
which coincides with the standard product metric on $S^3\times\R^4$ up to rescaling the first factor. The contact metric structure is given by
\[\alpha = -\Re p^{-1}dq\,j,\quad F=\frac12\Re \overline{p^{-1}dp}\wedge p^{-1}dq\,j=
\frac12\Re d\overline{p}\wedge dq\,j,\]
and the reduction to $\SU(3)$ by
\[\Omega=\frac{1}{54}(\overline pdp)^3 +\frac14\re \overline p dp (\overline p dq\,j)^2 + i\left(\frac18\re(\overline pdp)^2\overline pdq\,j +\frac16 (\overline p dq\,j)^3\right).\]
Now introduce the real coordinates
\[p=p_0+ip_1+jp_2+kp_3,\quad q=q_0+iq_1+jq_2+kq_3,\]
and consider the closed one-dimensional subgroup
\[K=\{(0,jq_2)\in \SU(2)\ltimes\h,q_2\in\R\}.\]
The left action of $K$ on $\SU(2)\ltimes\h$ preserves the $\SU(3)$-structure, and has the right-invariant vector field \[X=\frac{\partial}{\partial q_2}\]
as its associated fundamental vector field.
It follows that the moment map is
\[\mu = \alpha(X)=-(p_0dq_2+p_1dq_3-p_2dq_0-p_3dq_1)(X)=-p_0,\]
and $t=1$. Then $\mu^{-1}(0)=S^2\times\h$, where $S^2$ is the two-sphere of imaginary unit quaternions, and the contact reduction is well defined and smooth; explicitly,
\[(\SU(2)\ltimes\h)//K\cong S^2\times\R^3.\]
Since $\h$ is normal in $\SU(2)\ltimes\h$, at each point $X$ lies in the span of $e_2,e_3,e_5,e_7$. Hence, $Q(X)\equiv 0$, and by Theorem~\ref{thm:HypoToHypo}, we obtain a hypo-contact structure on $S^2\times\R^3$ with
\[B=\diag(0,-2,0,-2,0).\]

We can repeat the construction starting with $M=\Z^4\backslash (\SU(2)\ltimes_\phi\h)$, where $\Z^4\subset\h$ consists of points with  coordinates in $\Z$. Since $K$  and $\Z^4$ commute a free circle action is induced, and we can apply Theorem~\ref{thm:HypoToHypo}. To identify the contact reduction, observe that the above description of $\mu$ remains valid. Hence, we obtain a diagram
\[\xymatrix{
&\mu^{-1}(0)\ar[r]\ar[ld]^{S^1}\ar[dd]^{T^4}&M\ar[dd]_{T^4} \\
M//S^1\ar[dr]^{T^3}\\&S^2\ar[r] &\SU(2)}\]
with each arrow corresponding to a trivial torus bundle; in particular, \[M//S^1\cong S^2\times T^3\;.\]

\smallskip
\textbf{Acknowledgements}. We would like to thank S.~Salamon for helpful suggestions.
\bibliographystyle{plain}
\bibliography{contact}

\small\noindent Dipartimento di Matematica e Applicazioni, Universit\`a di Milano Bicocca,  Via Cozzi 53, 20125 Milano, Italy.\\
\texttt{diego.conti@unimib.it}

\vspace{0.15cm}
\small\noindent Dipartimento di Matematica, Universit\`a di Torino, Via Carlo Alberto 10, 10123 Torino, Italy.\\
\texttt{annamaria.fino@unito.it}
\end{document}